\documentclass[12pt]{article}
\usepackage{amssymb, amsmath, amsthm, amsfonts, tikz, geometry, ytableau, mathrsfs}
\usepackage[indent,margin=1cm]{caption}
\usepackage{fullpage}
\usepackage[colorlinks,linkcolor=blue,anchorcolor=blue,citecolor=blue]{hyperref}

\numberwithin{equation}{section}
\numberwithin{figure}{section}

\hypersetup{colorlinks = true}
\newtheorem{thm}{Theorem}[section]
\newtheorem{conj}[thm]{Conjecture}

\newtheorem{defi}[thm]{Definition}

\newtheorem{lem}[thm]{Lemma}
\newtheorem{prop}[thm]{Proposition}
\newtheorem{prob}[thm]{Problem}

\allowdisplaybreaks

\begin{document}
\begin{center}
	{\large \bf Strongly nice property and Schur positivity of graphs}
\end{center}

\begin{center}
Ethan Y.H. Li$^{1}$, Grace M.X. Li$^{2}$,  Arthur L.B. Yang$^{3}$ and Zhong-Xue Zhang$^{4}$\\[6pt]
\end{center}

\begin{center}
$^{1}$School of Mathematics and Statistics,\\
Shaanxi Normal University, Xi'an, Shaanxi 710119, P. R. China

$^{2}$School of Mathematics and Data Science, Shaanxi University of Science and
Technology, Xi’an, Shaanxi 710021, P. R. China

$^{3,4}$Center for Combinatorics, LPMC\\
Nankai University, Tianjin 300071, P. R. China\\[6pt]

Email: $^{1}${\tt yinhao\_li@snnu.edu.cn}, $^{2}${\tt grace\_li@sust.edu.cn}, $^3${\tt yang@nankai.edu.cn}, $^{4}${\tt zhzhx@mail.nankai.edu.cn}
\end{center}

\noindent\textbf{Abstract.}
Motivated by the notion of nice graphs, we introduce the concept of strongly nice property, which can be used to study the Schur positivity of symmetric functions. We show that a graph and all its induced subgraphs are strongly nice if and only if it is claw-free, which strengthens a result of Stanley and provides further evidence for the well-known conjecture on the Schur positivity of claw-free graphs.
As another application, we solve Wang and Wang's conjecture on the non-Schur positivity of squid graphs
$Sq(2n-1;1^n)$ for $n \ge 3$ by proving that these graphs are not strongly nice. \\[6pt]
\noindent \emph{AMS Mathematics Subject Classification 2020:} 05E05, 06A07

\noindent \emph{Keywords:} nice, strongly nice, Schur positivity, claw-free graphs, squid graphs

\section{Introduction}\label{intro}

Schur positivity is of great importance in combinatorics, since it has a deep relationship with representation theory and algebraic geometry. There have been plenty of conjectures on Schur positivity of certain symmetric functions. One of the most interesting conjectures is as follows, which  was firstly proposed by Gasharov (unpublished) and explicitly stated by Stanley \cite{Sta98}.
\begin{conj}[{\cite[Conjecture 1.4]{Sta98}}]\label{conj-clawfree}
The chromatic symmetric functions of all claw-free graphs (containing no induced subgraph isomorphic to the claw $K_{1,3}$) are Schur positive.
\end{conj}
However, in most cases it is very difficult to determine whether a symmetric function is Schur positive or not, which leads to further study on sufficient or necessary conditions for Schur positivity.
The nice property, which was defined for graphs by Stanley \cite{Sta98}, serves as a useful necessary condition for the Schur positivity of its chromatic symmetric function.
In particular, Stanley proved the following result.
%
\begin{prop}[{\cite[Proposition 1.5 and Proposition 1.6]{Sta98}}]\label{prop-nice}
If the chromatic symmetric function of a graph $G$ is Schur positive, then $G$ is nice. A graph $G$ is claw-free if and only if $G$ and all its induced subgraphs are nice.
\end{prop}

Proposition \ref{prop-nice} has the following interesting applications. On the one hand, as noted by Stanley \cite{Sta98}, it provides evidence for Conjecture \ref{conj-clawfree}. On the other hand,
one can use the above result to prove the non-Schur positivity of certain chromatic symmetric functions by showing that they are not nice.
For example, Dahlberg, She and van Willigenburg \cite{DSvW20} used this basic idea to prove that any $n$-vertex bipartite graph with a vertex of degree more than $\lceil \frac{n}{2} \rceil$ is not Schur positive. The same idea was also used by Wang and Wang \cite{WW20} to prove the non-Schur positivity of wheel graphs $W_n$ ($n\geq 7$), windmill graphs $W_{n}^d$ ($n,d\geq 3$), and complete bipartite graphs $K_{m,n}$  ($m\geq 4$). Li, Qiu, Yang and Zhang \cite{LQYZ24} constructed a family of distributive lattices being not nice, and thus answered an open problem proposed by Stanley \cite{Sta98}.

\begin{figure}[htbp]
\centering
\begin{tikzpicture}[scale = 1]

\fill (1,2) node (v1) {} circle (0.5ex);
\fill (1,0) node (v1) {} circle (0.5ex);
\fill (0,1) node (v1) {} circle (0.5ex);
\fill (-1,2) node (v1) {} circle (0.5ex);
\fill (-1,0) node (v1) {} circle (0.5ex);
\draw [thick] (-1.5,0.134) arc (-120.0007:40:1);
\draw [thick] (-0.234,1.6428) arc (40.0022:140:1);
\draw [thick] (-1.5,0.134) arc (-120.0007:-140:1);
\node at (-1.75,1.125) {$\vdots$};
\draw [thick] (1,0) -- (0,1) -- (1,2);
\node at (1,1.125) {$\vdots$};

\node [above] at (-1,2) {$u_1$};
\node [left] at (0,1) {$u$};
\node [below] at (-1,0) {$u_{2n-2}$};
\node [right] at (1,2) {$v_1$};
\node [right] at (1,0) {$v_n$};
\end{tikzpicture}
    \caption{The squid graph $Sq(2n-1;1^n)$}\label{fig-squid}
\end{figure}
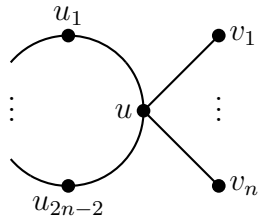

However, sometimes this approach does not work since there do exist many graphs which are nice but not Schur positive. In \cite{WW20} Wang and Wang studied the $s$-positivity of a class of squid graphs $Sq(2n-1; 1^n)$ defined by attaching $n$ leaves to one vertex of a cycle $C_{2n-1}$, as shown in Figure \ref{fig-squid}. They proposed the following conjecture.

\begin{conj}[{\cite[Conjecture 3.4]{WW20}}]\label{conj-squid}
The squid graph $Sq(2n-1; 1^n)$ is not Schur positive for $n \ge 3$.
\end{conj}

As will be shown in Theorem \ref{thm-squid-nice}, the graph $Sq(2n-1; 1^n)$ is nice for all $n \ge 3$. This implies that it is impossible to prove the non-Schur positivity of these graphs by showing they are not nice.

Motivated by Conjecture \ref{conj-clawfree}, Proposition \ref{prop-nice}, and Conjecture \ref{conj-squid}, we introduce the notion of strongly nice property for graphs and symmetric functions.
In Section \ref{sec-snice} we give the formal definition of strongly nice property, and show that Schur positivity implies strongly nice property.
In Section \ref{sec-claw} we strengthen Proposition \ref{prop-nice} to the strongly nice property. In Section \ref{sec-squid} we prove Conjecture \ref{conj-squid} by showing that $Sq(2n-1;1^n)$ is not strongly nice.
In Section \ref{sec-prob} we propose one question on the strongly nice property of incomparability graphs of Boolean lattices.

\section{Strongly nice property}\label{sec-snice}

This section is devoted to defining the strongly nice property for graphs and symmetric functions, and establishing its connection with Schur positivity. Recall that Stanley \cite{Sta98} defined the nice property first for graphs, and then for posets
by using their incomparability graphs. In this paper, we shall define the strongly nice property first for symmetric functions, and then for graphs and posets by their chromatic symmetric functions.

Let us begin with some basic definitions on symmetric functions. For more details, see \cite{Mac15} or \cite{StaEC2}.
Given a set of countably infinite indeterminates $\mathbf{x} = \{x_1,x_2,\ldots\}$, the algebra $\mathbb{Q}[[\mathbf{x}]]$ is defined to be the commutative algebra of formal power series in these indeterminates over the rational field $\mathbb{Q}$. The algebra of symmetric functions $\Lambda_{\mathbb{Q}}(\mathbf{x})$ is defined as the subalgebra of $\mathbb{Q}[[\mathbf{x}]]$ consisting of formal power series $f$ of bounded degree and satisfying
\[
f(\mathbf{x})=f(x_1,x_2,\ldots) = f(x_{\omega(1)},x_{\omega(2)},\ldots)
\]
for any permutation $\omega$ of positive integers. We usually abbreviate $f(\mathbf{x})$ to $f$ throughout this paper.

The bases of $\Lambda_{\mathbb{Q}}(\mathbf{x})$ are indexed by (integer) partitions. A partition of $n$ is a sequence $\lambda = (\lambda_1,\ldots,\lambda_\ell)$ satisfying
\[
\lambda_1 \ge \lambda_2 \ge \cdots \ge \lambda_\ell > 0 \quad \mbox{and} \quad \lambda_1 + \lambda_2 + \cdots + \lambda_\ell = n,
\]
where $\ell = \ell(\lambda)$ denotes the length of $\lambda$. By convention, we set $\lambda_i=0$ for $i>\ell(\lambda)$.
Given two partitions $\lambda,\mu$ of the same number $n$, we say $\lambda \ge \mu$ in dominance order if
$$\sum_{i=1}^k \lambda_i \ge \sum_{i=1}^k \mu_i$$
holds for all $k \ge 1$.

This paper is mainly concerned with two bases of $\Lambda_{\mathbb{Q}}(\mathbf{x})$: the monomial symmetric functions $m_{\lambda}$ and the Schur functions $s_{\lambda}$. For any partition $\lambda = (\lambda_1,\ldots,\lambda_\ell)$, the monomial symmetric function $m_{\lambda}$ is defined as
$$m_{\lambda}=\sum_{\alpha} x^{\alpha},$$
where $\alpha$ ranges over all distinct permutations of $(\lambda_1,\ldots,\lambda_\ell,0,0,\ldots)$ and $x^{\alpha} = x_1^{\alpha_1}x_2^{\alpha_2}\cdots$ for $\alpha = (\alpha_1,\alpha_2,\ldots)$. The Schur function $s_{\lambda}$ is defined as
\[
s_{\lambda} = \sum_{\mu} K_{\lambda\mu}m_{\mu},
\]
where $K_{\lambda\mu}$ denotes the number of semi-standard Young tableaux of shape $\lambda$ and type $\mu$. Given a symmetric function $f$ and a basis $\{b_{\lambda}\}$ of $\Lambda_{\mathbb{Q}}(\mathbf{x})$, let $[b_{\lambda}]f$ denote the coefficient of $b_{\lambda}$ in $f$.
A symmetric function $f$ is said to be \emph{Schur positive} or \emph{$s$-positive} if $[s_{\lambda}]f\geq 0$ for any partition $\lambda$.

The chromatic symmetric functions of graphs have been extensively studied since they were introduced by
Stanley \cite{Sta95}.
Let $G$ be a graph with vertex set $V(G) = \{v_1,\ldots,v_d\}$. Then the \textit{chromatic symmetric function} of $G$ is defined by Stanley as
\[
X_G = \sum_{\kappa} x_{\kappa(v_1)}\cdots x_{\kappa(v_d)},
\]
where $\kappa: V(G) \to \{1,2,\ldots\}$ ranges over all proper colorings of $G$, i.e., $\kappa(u) \neq \kappa(v)$ for any edge $uv \in E(G)$. As for any poset $P$, the chromatic symmetric function is defined on its \textit{incomparability graph} $\mathrm{inc}(P)$, whose vertex set consists of elements of $P$ and edge set is formed by pairs of vertices not comparable in $P$.

Stanley obtained a combinatorial expansion of $X_G$ in terms of monomial symmetric functions by using stable partitions of $V(G)$. By a \emph{stable partition} of $G$ we mean a set partition $B = \{B_1,\ldots,B_k\}$ of $V(G)$ such that any pair of vertices in the same block $B_i$ ($1 \le i \le k$) are not adjacent (or equivalently, $B_i$ is a stable set). A \textit{semi-ordered stable partition} is obtained by ordering the blocks of the same size. For instance, taking $G$ to the empty graph on five vertex set  $V(G)=\{1,2,3,4,5\}$, one should consider
$\{\{2,3\},\{4,5\},\{1\}\}$ and $\{\{4,5\},\{2,3\},\{1\}\}$ as the same stable partition of $V(G)$ but as two different semi-ordered stable partitions. The type of a (semi-ordered) stable partition $B$ is defined to be the integer partition formed by rearranging the block sizes $|B_1|,\ldots,|B_k|$ in weakly decreasing order. Stanley obtained the following result.

\begin{prop}\cite[Proposition 2.4]{Sta95}\label{prop-mexp}
Let $\tilde{a}_{\lambda}$ be the number of semi-ordered stable partitions of $G$ of type $\lambda$. Then
\[
X_G = \sum_{\lambda} \tilde{a}_{\lambda} m_{\lambda}.
\]
\end{prop}

Stanley \cite{Sta98} showed that the Schur positivity of $X_G$ can be used to study the nice property of $G$, and vice versa. A graph $G$ is called \textit{nice} if, for any pair of partitions $\lambda,\mu$ satisfying $\lambda \ge \mu$ in dominance order, the graph $G$ must contain a stable partition of type $\mu$ as long as $G$ contains a stable partition of type $\lambda$.
By Proposition \ref{prop-mexp}, a graph $G$ is  {nice}
if and only if whenever the coefficient $\tilde{a}_{\lambda}$ of $m_{\lambda}$ in $X_G$ does not vanish and whenever $\lambda \ge \mu$ in dominance order, then the coefficient $\tilde{a}_{\mu}$ does not vanish.
In this manner, the nice property can be naturally defined for any symmetric function. Precisely, we say that a symmetric function $f$ is \textit{nice} if for any pair of partitions $\mu \le \lambda$  in dominance order with $[m_{\lambda}]f > 0$ we have $[m_{\mu}]f > 0$.

In the following we strengthen the nice property of symmetric functions to a quantitive version.

\begin{defi}
A symmetric function $f$ is said to be \textup{strongly nice} if $[m_{\mu}]f \ge [m_{\lambda}]f$ whenever $\mu \le \lambda$ in dominance order.
\end{defi}

We say that a graph $G$ is  strongly nice if $X_G$ is strongly nice, or equivalently, if the number of semi-ordered stable partitions of $G$ of type $\mu$ is more than or equal to that of type $\lambda$ for any pair $\mu \le \lambda$. The following result implies that strongly nice property is more powerful than nice property for studying Schur positivity or non-Schur positivity.

\begin{lem}\label{lem-snice}
A strongly nice symmetric function is always nice, and an $s$-positive symmetric function is always strongly nice.
\end{lem}

\begin{proof}
The first assertion follows directly from the definitions. To prove the second, we need to use the monotonicity of the Kostka numbers due to White \cite{Whi80}, who proved that $K_{\nu\mu} \ge K_{\nu\lambda}$ whenever $\mu \le \lambda$ in dominance order. Suppose that $f=\sum_{\nu}c_{\nu}s_{\nu}$ with $c_{\nu} \ge 0$ for all $\nu$.
Then
\[
[m_{\mu}]f = [m_{\mu}]\sum_{\nu}c_{\nu}s_{\nu} = [m_{\mu}]\sum_{\nu}c_{\nu}\sum_{\rho}K_{\nu\rho}m_{\rho} = \sum_{\nu}c_{\nu}K_{\nu\mu} \ge \sum_{\nu}c_{\nu}K_{\nu\lambda} = [m_{\lambda}]f.
\]
This completes the proof.
\end{proof}

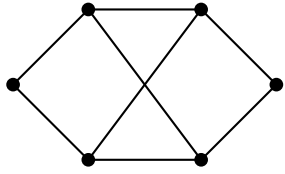
\begin{figure}[htbp]
\centering
\begin{tikzpicture}[scale = 1]
\fill (0,1) node (v1) {} circle (0.5ex);
\fill (-1,0) node (v9) {} circle (0.5ex);
\fill (1.5,1) node (v10) {} circle (0.5ex);
\fill (2.5,0) node (v10) {} circle (0.5ex);
\fill (1.5,-1) node (v10) {} circle (0.5ex);

\draw [thick] (0,1) -- (-1,0);
\draw [thick] (0,1) -- (1.5,1) -- (2.5,0) -- (1.5,-1);
\draw [thick] (-1,0) -- (0,-1);
\draw [thick] (1.5,1) -- (0,-1);
\draw [thick] (1.5,-1) -- (0,1);
\fill (0,-1) node (v10) {} circle (0.5ex);
\draw [thick] (0,-1) -- (1.5,-1);

\end{tikzpicture}
    \caption{A strongly nice graph without Schur positivity}\label{fig-nice-nonspos}
\end{figure}

\noindent\textbf{Remark.} Similar to the nice property, the strongly nice property is also not equivalent to $s$-positivity, even in the special case of chromatic symmetric functions. Figure \ref{fig-nice-nonspos} presents a graph which is strongly nice but not $s$-positive, whose chromatic symmetric functions is calculated by SageMath \cite{Sage} as
{\small \begin{align*}
&720m_{(1,1,1,1,1,1)}+ 168m_{(2,1,1,1,1)}+ 44m_{(2,2,1,1)}+ 6m_{(2,2,2)}+ 12m_{(3,1,1,1)}+6m_{(3,2,1)}+ 2m_{(3,3)}\\
=&152s_{(1,1,1,1,1,1)} + 52s_{(2,1,1,1,1)} + 26s_{(2,2,1,1)} - 4s_{(2,2,2)} + 2s_{(3,1,1,1)} + 4s_{(3,2,1)} + 2s_{(3,3)}.
\end{align*}}

\section{Claw-free graphs}\label{sec-claw}

In this section we shall study the strongly nice property of claw-free graphs.
The main result of this section is as follows, which strengthens Proposition \ref{prop-nice} and provides further evidence for Conjecture \ref{conj-clawfree}.

\begin{thm}\label{thm-strongly-nice}
A graph $G$ is claw-free if and only if $G$ and all its induced subgraphs are strongly nice.
\end{thm}

\begin{proof}
The sufficiency is straightforward since the claw graph is not strongly nice. Indeed, the chromatic symmetric function of the claw graph is
\[
X_{K_{1,3}} = 24m_{(1,1,1,1)} + 6m_{(2,1,1)} + m_{(3,1)},
\]
while $(2,2)<(3,1)$ but $[m_{(2,2)}]X_{K_{1,3}} = 0 < 1 = [m_{(3,1)}]X_{K_{1,3}}$.

Now we proceed to prove the necessity. It suffices to prove the following claim since any induced subgraph of a claw-free graph is also claw-free.

\textbf{Claim.} If a graph $G$ is claw-free, then it is strongly nice.

To prove the strongly nice property of claw-free graphs, by definition, we only need to show that $[m_{\mu}]X_G \ge [m_{\lambda}]X_G$ for all partitions $\mu,\lambda$ with $\lambda$ covering $\mu$ under dominance order. Such partitions are characterized as follows: if $\mu$ is covered by $\lambda$, then there exists $i < j$ such that $\mu_i = \lambda_i-1$, $\mu_j = \lambda_j + 1$, and $\mu_k = \lambda_k$ for $k \neq i,j$. Note that $0\leq \lambda_j \le \lambda_i - 2$ since $\mu_i \ge \mu_j \ge 1$.
Clearly, if $\mu$ is covered by $\lambda$, then $\ell(\mu)=\ell(\lambda)$ or
$\ell(\mu)=\ell(\lambda)+1$.

If $[m_{\lambda}]X_G = 0$, then the inequality naturally holds. From now on, we assume that $[m_{\lambda}]X_G > 0$.
By Proposition \ref{prop-mexp}, $[m_{\lambda}]X_G$ is equal to $\tilde{a}_{\lambda}$, the cardinality of the set $\tilde{A}_{\lambda}$ of semi-ordered stable partitions of type $\lambda$. Hence in the following we shall establish an injection $\phi$ from $\tilde{A}_{\lambda}$ to $\tilde{A}_{\mu}$.
For notational convenience, we use the unique representation of semi-ordered stable partitions $B = \{B_1,\ldots,B_{\ell(\lambda)}\}$ of type $\lambda$, which is obtained by arranging the blocks of difference sizes in weakly decreasing order, or equivalently, requiring that $|B_k|=\lambda_k$ for all $1 \le k \le \ell(\lambda)$.


Consider the subgraph $H$ induced by $B_i \cup B_j$,
where $i,j$ are uniquely determined by $\lambda$ and $\mu$ as mentioned above.
If $\lambda_j=0$, then we set $B_j=\emptyset$ for convenience, though it is no longer a valid block of the semi-ordered stable partition $B$.
It is clear that $H$ is a claw-free bipartite graph. Then the maximum degree of $H$ is at most 2 and $H$ must be a disjoint union of paths and (even) cycles, where we regard isolated vertices as paths. It follows from $\lambda_i > \lambda_j$ that there exists at least one odd path $P$ (containing an odd number of vertices) with $|P\cap B_i| > |P\cap B_j|$ (actually $|P\cap B_i| = |P\cap B_j|+1$). 
Now fix an arbitrary labeling $\alpha$ of $V(G)$. For each odd path $P$ in $H$, define
\[
\alpha(P) = \min \{\alpha(v) \mid v \in P\}.
\]
We may assume that such odd paths are $P_1,\ldots,P_t$ with $\alpha(P_1)<\cdots<\alpha(P_t)$.
Then define the word of $B$ as
\[
W(B) = c_1c_2 \cdots c_t, \quad \mbox{where}\ c_k = \begin{cases}
                                 1, & \mbox{if } |P_{k} \cap B_i| > |P_{k} \cap B_j| \\
                                 2, & \mbox{if } |P_{k} \cap B_i| < |P_{k} \cap B_j|
                               \end{cases} \quad\mbox{ for } 1 \le k \le t.
\]
For example, the semi-ordered stable partition shown in Figure \ref{fig-word} (we only present $B_i$ and $B_j$ for convenience) has word $W(B) = 1211$.
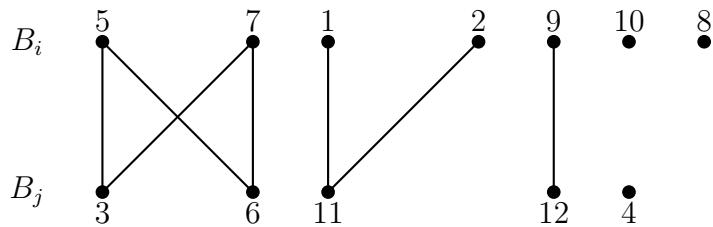
\begin{figure}[htbp]
\centering
\begin{tikzpicture}[scale = 1]

\fill (0,1) node (v1) {} circle (0.5ex);
\fill (0,-1) node (v9) {} circle (0.5ex);
\fill (2,-1) node (v9) {} circle (0.5ex);
\fill (2,1) node (v10) {} circle (0.5ex);

\draw [thick] (2,-1) -- (0,1) -- (0,-1);
\draw [thick] (0,-1) -- (2,1) -- (2,-1);

\fill (3,-1) node (v9) {} circle (0.5ex);
\fill (3,1) node (v10) {} circle (0.5ex);
\fill (5,1) node (v10) {} circle (0.5ex);

\draw [thick] (3,1) -- (3,-1) -- (5,1);

\fill (6,-1) node (v9) {} circle (0.5ex);
\fill (6,1) node (v10) {} circle (0.5ex);

\draw [thick] (6,1) -- (6,-1);

\fill (7,1) node (v9) {} circle (0.5ex);
\fill (7,-1) node (v9) {} circle (0.5ex);
\fill (8,1) node (v10) {} circle (0.5ex);

\node [above] at (0,1) {5};
\node [above] at (2,1) {7};
\node [below] at (0,-1) {3};
\node [below] at (2,-1) {6};

\node [above] at (3,1) {1};
\node [above] at (5,1) {2};
\node [above] at (6,1) {9};
\node [above] at (7,1) {10};
\node [above] at (8,1) {8};
\node [below] at (7,-1) {4};

\node [below] at (3,-1) {11};
\node [below] at (6,-1) {12};
\node at (-1,1) {$B_i$};
\node at (-1,-1) {$B_j$};
\end{tikzpicture}
    \caption{The blocks $B_i$ and $B_j$}\label{fig-word}
\end{figure}
\begin{figure}[htbp]
\centering
\begin{tikzpicture}[scale = 1]

\fill (0,1) node (v1) {} circle (0.5ex);
\fill (0,-1) node (v9) {} circle (0.5ex);
\fill (2,-1) node (v9) {} circle (0.5ex);
\fill (2,1) node (v10) {} circle (0.5ex);

\draw [thick] (2,-1) -- (0,1) -- (0,-1);
\draw [thick] (0,-1) -- (2,1) -- (2,-1);

\fill (3,-1) node (v9) {} circle (0.5ex);
\fill (3,1) node (v10) {} circle (0.5ex);
\fill (5,1) node (v10) {} circle (0.5ex);

\draw [thick] (3,1) -- (3,-1) -- (5,1);

\fill (6,-1) node (v9) {} circle (0.5ex);
\fill (6,1) node (v10) {} circle (0.5ex);

\draw [thick] (6,1) -- (6,-1);

\fill (7,1) node (v9) {} circle (0.5ex);
\fill (7,-1) node (v9) {} circle (0.5ex);
\fill (8,-1) node (v10) {} circle (0.5ex);

\node [above] at (0,1) {5};
\node [above] at (2,1) {7};
\node [below] at (0,-1) {3};
\node [below] at (2,-1) {6};

\node [above] at (3,1) {1};
\node [above] at (5,1) {2};
\node [above] at (6,1) {9};
\node [above] at (7,1) {8};
\node [below] at (8,-1) {10};
\node [below] at (7,-1) {4};

\node [below] at (3,-1) {11};
\node [below] at (6,-1) {12};
\node at (-1,1) {$\bar{B}_i$};
\node at (-1,-1) {$\bar{B}_j$};
\end{tikzpicture}
    \caption{The blocks $\bar{B}_i$ and $\bar{B}_j$}\label{fig-word-2}
\end{figure}
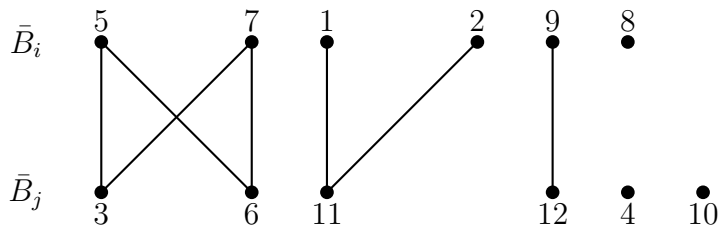
The desired map $\phi$ is constructed as follows. Let $w_1(W)$ and $w_2(W)$ be the number of 1 and 2 of a word $W$, respectively. Let $1 \le p \le t$ be the smallest index such that $w_1(c_1\cdots c_p) - w_2(c_1\cdots c_p)$ is maximum. Then we have $w_1(c_1\cdots c_p) - w_2(c_1\cdots c_p) \ge w_1(c_1\cdots c_t) - w_2(c_1\cdots c_t) \ge 2$ since $|B_i| = \lambda_i \ge \lambda_j+2 = |B_j|+2$, which implies $p \ge 2$. If follows that $c_p=1$ since otherwise the desired index would be $p-1$. Define $\phi(B)=\bar{B}$, where $\bar{B}$ is the semi-ordered stable partition of $V(G)$ obtained from $B$ by exchanging the vertices of $P_p$ in $B_i$ and $B_j$, or precisely, by letting
$$\bar{B}_i=
(B_i \setminus V(P_{p})) \cup (B_j \cap V(P_{p})),\quad
\bar{B}_j=(B_j \setminus V(P_{p})) \cup (B_i \cap V(P_{p}))$$
and fixing the remaining blocks. One can verify that $\bar{B} \in \tilde{A}_{\mu}$. If we define $W(\bar{B})=\bar{c}_1\bar{c}_2 \cdots \bar{c}_t$ in the same way as $W({B})$, it is clear that $\bar{c}_p=2$ and $\bar{c}_l=c_l$ for all $l\neq p$. Moreover, we have
\begin{align}\label{eq-ccbar}
w_1(\bar{c}_1\cdots \bar{c}_k) - w_2(\bar{c}_1\cdots \bar{c}_k)
=\begin{cases}
w_1({c}_1\cdots {c}_k) - w_2({c}_1\cdots {c}_k), & \mbox{if } k<p; \\
w_1({c}_1\cdots {c}_k) - w_2({c}_1\cdots {c}_k)-2, & \mbox{if } k\geq p.
 \end{cases}
\end{align}

We proceed to show that $\phi$ is injective.
To this end, we construct a map $\varphi$ on the image set $\phi(\tilde{A}_{\lambda})$ and prove that the composition $\varphi\circ\phi$ is the identity map on $\tilde{A}_{\lambda}$.
Note that $\tilde{A}_{\lambda}$ can be divided into disjoint subsets according to the set $B_i \cup B_j$, namely, $\tilde{A}_{\lambda} = \biguplus_H \tilde{A}_{\lambda}^H$, where $H$ ranges over all induced subgraphs in
$$ \{G[B_i \cup B_j] \mid B_i,B_j \in B \mbox{ for some semi-ordered stable partion $B$} \},$$
and $\tilde{A}_{\lambda}^H$ denotes the set of semi-ordered stable partition $B$ of type $\lambda$ with $B_i \cup B_j = V(H)$. One can further observe that $\phi(\tilde{A}_{\lambda}^H) \subseteq \tilde{A}_{\mu}^H$ since $B_i \cup B_j = \bar{B}_i \cup \bar{B}_j$. Hence it suffices to prove that the restriction of $\phi$ on $\tilde{A}_{\lambda}^H$ is injective for any $H$.

When fixing $H$, the labeling and the definition for the word remain the same.
Given a semi-ordered stable partition $\bar{B} = \{\bar{B}_1,\ldots,\bar{B}_{\ell(\mu)}\}$ in
$\phi(\tilde{A}_{\lambda}^H)$, consider the word
$W(\bar{B})=\bar{c}_1\bar{c}_2 \cdots \bar{c}_t$ and
choose the largest index $q$ such that the number $w_1(\bar{c}_1\cdots \bar{c}_q) - w_2(\bar{c}_1\cdots \bar{c}_q)$ is maximum. By \eqref{eq-ccbar} we have $q \le t-1$ and $\bar{c}_{q+1}=2$.
Then define $\varphi(\bar{B})=\hat{B}$, where $\hat{B}$ is the semi-ordered stable partition of $V(G)$ obtained from $\bar{B}$ by letting
$$\hat{B}_i=
(\bar{B}_i \setminus V(P_{q+1})) \cup (\bar{B}_j \cap V(P_{q+1})),\quad
\hat{B}_j=(\bar{B}_j \setminus V(P_{q+1})) \cup (\bar{B}_i \cap V(P_{q+1}))$$
and fixing the remaining blocks. Similarly, define $W(\hat{B})=\hat{c}_1\hat{c}_2 \cdots \hat{c}_t$. One can verify that $\hat{c}_{q+1}=1$ and $\hat{c}_l=\bar{c}_l$ whenever $l\neq q+1$. For example, the word in Figure \ref{fig-word-2} is $W(\bar{B}) = 1212$ and $W(\hat{B})) = 1211$.  We would like to mention that if $\mu_j=1$ then $\hat{B}_j$ is empty.
%
By the construction of $\phi$ and $\varphi$, one can check that
$p=q+1$ if $\bar{B} = \phi(B)$. It follows that $\varphi(\phi(B)) = \varphi(\bar{B}) = B$, implying the injectivity of $\phi$ on $\tilde{A}_{\lambda}^H$. This completes the proof.
\end{proof}

\noindent\textbf{Remark.} In the proof of Proposition \ref{prop-nice}, Stanley only treated stable partitions since the main focus is the existence and the order is irrelevant. However, in our proof, we have to make use of the semi-order to define the desired injection, and different blocks of the same size are treated differently according to the order.

\section{Squid graphs}\label{sec-squid}

The main objective of this section is to prove Conjecture \ref{conj-squid}. Precisely, we have the following result.


\begin{thm}\label{thm-squid}
For $n \ge 3$, the squid graph $Sq(2n-1;1^n)$ is not strongly nice. Moreover, $Sq(2n-1;1^n)$ is not $s$-positive.
\end{thm}

\begin{proof}
Precisely, we are going to show
\begin{equation}\label{eq-sq-ineq}
[m_{(n,n,n-1)}]X_{Sq(2n-1;1^n)} < [m_{(n+1,n-1,n-1)}]X_{Sq(2n-1;1^n)}.
\end{equation}
Using the labeling in Figure \ref{fig-squid}, the key observation is that in any stable partition the size of the block containing $u$ is at most $n-1$. Indeed, the vertices $v_1,\ldots,v_n$ cannot be in a block containing $u$, and there are at most $n-2$ vertices in $\{u_2,\ldots,u_{2n-3}\}$ which can be put into this block. Moreover, it is not difficult to check that there are exactly $n-1$ ways to choose such a block of size $n-1$ (choosing $n-2$ non-adjacent points in $\{u_2,\ldots,u_{2n-3}\}$).

Once we have selected a stable set of size $n-1$ containing $u$, deleting these vertices will result in a subgraph $L$ consisting of one edge and $2n-2$ isolated vertices. Therefore, the number of ways to obtain a semi-ordered stable partition of type $(n,n)$ in $L$ is $\binom{2}{1} \cdot \binom{2n-2}{n-1} = 2\binom{2n-2}{n-1}$. Similarly, the number of ways to obtain a semi-ordered stable partition of type $(n+1,n-1)$ in $L$ is $\binom{2}{1} \cdot \binom{2n-2}{n} = 2\binom{2n-2}{n}$. Hence,
\begin{align*}
[m_{(n,n,n-1)}]X_{Sq(2n-1;1^n)} &= 2(n-1)\binom{2n-2}{n-1},\\[5pt]
[m_{(n+1,n-1,n-1)}]X_{Sq(2n-1;1^n)} &= 4(n-1)\binom{2n-2}{n},
\end{align*}
where the second equality is obtained by distinguishing the two blocks of size $n-1$.
Now \eqref{eq-sq-ineq} follows from
\[
\frac{2(n-1)\binom{2n-2}{n-1}}{4(n-1)\binom{2n-2}{n}} = \frac{n}{2(n-1)} < 1,
\]
which is valid exactly for $n \ge 3$.
\end{proof}


\begin{thm}\label{thm-squid-nice}
The squid graph $Sq(2n-1; 1^n)$ is nice for $n \ge 3$.
\end{thm}
\begin{proof}
Label the vertices of $Sq(2n-1; 1^n)$ as in Figure \ref{fig-squid}.  Observe that $$\{\{u_1,\ldots,u_{2n-3},v_1,\ldots,v_n\},\{u_2,\ldots,u_{2n-2}\},\{u\}\}$$ is a stable partition of $Sq(2n-1; 1^n)$ of type $\lambda = (2n-1,n-1,1)$. Let $\mu=(\mu_1,\ldots,\mu_\ell)$ be any partition of $3n-1$.
We claim that there exists a stable partition of $Sq(2n-1; 1^n)$ of type $\mu$ if and only if $\mu\le\lambda$ in dominance order, which would imply its nice property by definition.

Suppose that there exists a stable partition of $Sq(2n-1; 1^n)$ of type $\mu$. Then any stable set of $Sq(2n-1; 1^n)$ has size at most $2n-1$ and hence $\mu_1\le 2n-1$. Note that any odd cycle cannot be divided into two stable sets since it is not bipartite, which yields $\mu_3 \ge 1$. It follows that $\mu_1+\mu_2\le 3n-1-1=3n-2=\lambda_1+\lambda_2$, and hence $\mu\le\lambda$.

Conversely, we assume that $\mu\le\lambda$. Then $\ell = \ell(\mu) \ge 3$ and $\mu_{\ell} \le n-1$ since otherwise $\mu_1 + \cdots + \mu_{\ell} \ge \ell n > 3n-1$. Now we proceed to construct a stable partition $B = \{B_1,B_2,\ldots,B_\ell\}$ of type $\mu$ with $|B_k|=\mu_k$ for each $1\leq k\leq \ell$. At first, we take $B_\ell=\{u,u_2,\ldots,u_{2(\mu_{\ell}-1)}\}$. If $\mu_1 > n-\mu_l$, then we set
\small \[
B_{1}=\begin{cases}
        \{u_{2\mu_\ell}, u_{2(\mu_\ell+1)},\ldots, u_{2n-2},v_1,v_2,\ldots,v_{\mu_1+\mu_\ell-n}\}, & \mbox{if } \mu_1 \le 2n-\mu_\ell,  \\
        \{u_{2\mu_\ell}, u_{2(\mu_\ell+1)},\ldots, u_{2n-2},v_1,v_2,\ldots,v_n,u_1,u_3,\ldots,u_{2(\mu_1+\mu_\ell-2n)-1}\}, & \mbox{if } \mu_1 > 2n-\mu_\ell,
      \end{cases}
\]
where the first is a valid stable set since in this case $\mu_1+\mu_\ell-n \le n$, and the second follows from $\mu_1 \le 2n-1$ and $2(\mu_1+\mu_\ell-2n)-1 \le 2(\mu_l-1)-1$. Then $V(Sq(2n-1;1^n)) \setminus (B_1 \cup B_\ell)$ consists of only isolated points, and hence the other blocks can be chosen arbitrarily.
If $\mu_1 \le n-\mu_\ell$, then there exists an index $i \ge 1$ such that $\mu_1 + \cdots + \mu_{i} \le n-\mu_\ell$ and $\mu_1 + \cdots + \mu_{i} + \mu_{i+1} > n-\mu_\ell$. Since $2(\mu_\ell+\mu_1+\cdots+\mu_{i}-1) \le 2(n-1)$, we can take
\[
B_j = \{u_{2(\mu_\ell+\mu_1+\cdots+\mu_{j-1})},u_{2(\mu_\ell+\mu_1+\cdots +\mu_{j-1}+1)},\ldots, u_{2(\mu_\ell+\mu_1+\cdots+\mu_{j}-1)}\}
\]
for $1 \le j \le i$. Then we set
\[
B_{i+1}=\{u_{2(\mu_\ell+\mu_1+\cdots+\mu_i)}, u_{2(\mu_\ell+\mu_1+\cdots+\mu_i)+2},\ldots, u_{2n-2},v_1,v_2,\ldots,v_{\mu_\ell + \mu_1+\cdots \mu_{i+1} -n}\},
\]
which is also possible since $\mu_{i+1} \le \mu_i \le n-\mu_\ell$ and
$$1\leq (\mu_\ell + \mu_1+\cdots +\mu_i) + \mu_{i+1} -n \le n + (n-\mu_\ell) -n \le n.$$
Now the graph induced by $V(Sq(2n-1;1^n)) \setminus (B_1 \cup \cdots \cup B_{i+1} \cup B_{\ell})$ consists of only isolated points, and hence the other blocks $B_{i+2},\ldots,B_{\ell-1}$ could be chosen arbitrarily. This completes the proof.
\end{proof}

Together with Theorem \ref{thm-squid}, the above result shows that squid graphs $Sq(2n-1; 1^n)$ form an infinite family of nice graphs which are not strongly nice.

\section{One open problem}\label{sec-prob}

Griggs \cite{Gri88} conjectured that the incomparability graph $\mathrm{inc}(B_n)$ of the Boolean lattice $B_n$ is nice. Stanley \cite{Sta98} further asked whether these graphs are $s$-positive. Stanley \cite{Sta98} noted the Schur positivity of $B_n$ for $n \le 4$, which implies that it is strongly nice. We already verified the nice property of $B_5$ by using SageMath \cite{Sage}. It is natural to ask the following problem.

\begin{prob}
Is $\mathrm{inc}(B_n)$ strongly nice?
\end{prob}

%
%
%
%

\noindent \textbf{Acknowledgments.}

Ethan Li is supported by the Fundamental Research Funds for the Central Universities (GK202207023). Arthur Yang is supported in part by the National Science Foundation of China (11971249 and 12325111).

\end{document}